\documentclass[12pt,reqno]{amsart}

\usepackage{amscd}
\usepackage{epsfig}
\usepackage{float}
\usepackage{url}
\usepackage[all]{xy}
\newcommand{\N}{{\mathbb N}}
\newcommand{\Z}{{\mathbb Z}}
\newcommand{\Q}{{\mathbb Q}}
\newcommand{\C}{{\mathbb C}}
\newcommand{\R}{{\mathbb R}}

\newcommand{\PP}{{\mathcal P}}

\newcommand{\www}{\widetilde}

\newcommand{\oooo}{\overline}

\newcommand{\paa}{\partial}

\newcommand{\Or}{{\rm Or}}

\DeclareMathOperator{\Aut}{Aut}
\DeclareMathOperator{\divis}{div}
\DeclareMathOperator{\lcm}{lcm}
\DeclareMathOperator{\ord}{ord}
\DeclareMathOperator{\supp}{supp}
\DeclareMathOperator{\tr}{tr}

\begin{document}

\theoremstyle{plain}
\newtheorem{lemma}{Lemma}[section]
\newtheorem{definition/lemma}[lemma]{Definition/Lemma}
\newtheorem{theorem}[lemma]{Theorem}
\newtheorem{proposition}[lemma]{Proposition}
\newtheorem{corollary}[lemma]{Corollary}
\newtheorem{conjecture}[lemma]{Conjecture}
\newtheorem{conjectures}[lemma]{Conjectures}

\theoremstyle{definition}
\newtheorem{definition}[lemma]{Definition}
\newtheorem{withouttitle}[lemma]{}
\newtheorem{remark}[lemma]{Remark}
\newtheorem{remarks}[lemma]{Remarks}
\newtheorem{example}[lemma]{Example}
\newtheorem{examples}[lemma]{Examples}
\newtheorem{notations}[lemma]{Notations}

\title[The combinatorics of weight systems]
{The combinatorics of weight systems and characteristic
polynomials of isolated quasihomogeneous singularities} 

\author{Claus Hertling and Makiko Mase}

\address{Claus Hertling\\
Lehrstuhl f\"ur algebraische Geometrie, 
Universit\"at Mannheim,
B6 26, 68159 Mannheim, Germany}

\email{hertling@math.uni-mannheim.de}

\address{Makiko Mase\\
Lehrstuhl f\"ur algebraische Geometrie, 
Universit\"at Mannheim,
B6 26, 68159 Mannheim, Germany}

\email{mmase@mail.uni-mannheim.de}

\date{August 02, 2021}

\subjclass[2010]{32S40, 05C25, 06A05, 11H56}

\keywords{quasihomogeneous singularity, weight system, 
monodromy, characteristic polynomial, combinatorial problems, 
Orlik blocks}

\thanks{This work was funded by the Deutsche 
Forschungsgemeinschaft (DFG, German Research Foundation) 
-- 242588615}


\begin{abstract}
{\Small 
A paper of the first author and Zilke proposed seven 
combinatorial problems around 
formulas for the characteristic polynomial and the exponents
of an isolated quasihomogeneous singularity. 
The most important of them was a conjecture on the 
characteristic polynomial. Here the conjecture is proved,
and some of the other problems are solved, too.  
In the cases where also an old conjecture of Orlik on the
integral monodromy holds, this has implications on the
automorphism group of the Milnor lattice.
The combinatorics used in the proof of the conjecture
consists of tuples of orders on sets $\{0,1,...,n\}$ with
special properties and may be of independent interest.}
\end{abstract}

\maketitle

\tableofcontents

\setcounter{section}{0}

\section{Introduction}\label{s1}
\setcounter{equation}{0}

\noindent
The paper \cite{HZ19} proposed seven combinatorial problems 
around formulas for the characteristic polynomial and the
exponents of an isolated quasihomogeneous singularity. 
Problem 6 was a conjecture on the characteristic polynomial. 
It was an amendment to an old conjecture of Orlik on the
integral monodromy (which is recalled below). 
Here the conjecture on the characteristic polynomial 
is proved, and also some  of the
other combinatorial problems are solved. 

The combinatorics which is used in the proof of the 
conjecture was found in \cite{HM20-2} and was used there
for a partial proof of Orlik's conjecture. 
It is described in detail in section \ref{s2}.
It may be of independent interest.

In this introduction the conjecture and some background
are explained. This contains problem 6 in \cite{HZ19}. 
Problem 7 in \cite{HZ19} is treated in section \ref{s3}.
Sections \ref{s4} sets some notations and notions.
Section \ref{s5}  solves problem 1. 
Section \ref{s6} discuss the combinatorics of the 
weight systems and of  the characteristic polynomial 
of an isolated quasihomogeneous singularity.
It solves problem 4 and discusses the problems 2, 3 and 5.

We start with the notion of an {\it Orlik block} and 
a result on its automorphisms.
Then we introduce isolated quasihomogenous singularities.
Finally we give Theorem \ref{t1.4},
which solves problem 6 in \cite{HZ19}, and Orlik's conjecture.

\begin{definition}\label{t1.1}
Let $M\subset\N=\{1,2,3,...\}$ be a finite nonempty subset. 
Its {\it Orlik block} is a pair $(H_M,h_M)$ with 
$H_M$ a $\Z$-lattice of rank $\sum_{m\in M}\varphi(m)$ and 
$h_M:H_M\to H_M$ an automorphism with characteristic polynomial
$\prod_{m\in M}\Phi_m$ ($\Phi_m$ is the $m$-th cyclotomic polynomial)
and with a cyclic generator $e_1\in M$, i.e.
\begin{eqnarray}\label{1.1}
H_M =\bigoplus_{j=1}^{rk M}\Z\cdot h_M^{j-1}(e_1).
\end{eqnarray}
$(H_M,h_M)$ is unique up to isomorphism.
$\Aut_{S^1}(H_M,h_M)$ denotes the group of all automorphisms of $H_M$
which commute with $h_M$ and which have all eigenvalues in $S^1$.
\end{definition}

\begin{definition}\label{t1.2}
A finite set $M\subset\N$ is enriched as follows to a 
directed graph $(M,E(M))$ with set of vertices $M$ 
and set of oriented edges $E(M)\subset M\times M$. 
An edge goes from $m_1\in M$ to $m_2\in M$ if 
$\frac{m_1}{m_2}$ is a power of a prime number $p$. 
Then it is called a $p$-edge.
\end{definition}

The main result in \cite{He20} is as follows.

\begin{theorem}\label{t1.3} \cite[Theorem 1.2]{He20}
Let $(H_M,h_M)$ be the Orlik block of a finite nonempty subset $M\subset\N$.
Then $\Aut_{S^1}(H_M,h_M)=\{\pm h_M^k\, |\, k\in\Z\}$ if and only if condition (I) or condition (II) in Definition 
\ref{t2.6} (c) and (d) are satisfied.
They are conditions on the graph $(M,E(M))$.
\end{theorem}

A weight system ${\bf w}=(w_1,...,w_n;1)$ with 
$w_i\in\Q_{>0}\cap (0,1)$
equips any monomial ${\bf x}^{\bf j}=x_1^{j_1}...x_n^{j_n}$ 
with a weighted degree 
$\deg_{\bf w}{\bf x}^{\bf j}:=\sum_{i=1}^n w_ij_i$. 
A polynomial $f\in\C[x_1,...,x_n]$ is an {\it isolated 
quasihomogeneous singularity} if for some weight system 
${\bf w}$ each monomial in $f$ 
has weighted degree 1 and if the functions 
$\frac{\paa f}{\paa x_1},...,\frac{\paa f}{\paa x_n}$ vanish 
simultaneously only at $0\in\C^n$. 
Then the {\it Milnor lattice} $H_{Mil}:=H_{n-1}(f^{-1}(1),\Z)$ 
(resp. the reduced homology in the case $n=1$) 
is a $\Z$-lattice of some rank $\mu\in\N$ \cite{Mi68}, 
which is called {\it Milnor number}. 
It comes equipped with a natural automorphism
$h_{Mil}:H_{Mil}\to H_{Mil}$ of finite order, 
the {\it monodromy}.
Its characteristic polynomial has the shape 
\begin{eqnarray*}
p_{ch,h_{Mil}} &= \prod_{m\in \N}\Phi_m^{\psi(m)}
\end{eqnarray*}
for a function $\psi:\N\to\N_0$ with finite support
$M_1:=\{n\in\N\,|\, \psi(n)\neq 0\}$. 
Denote $l_\psi:=\max(\psi(m)\, |\, m\in\N)$ and for 
$j=1,...,l_\psi$
\begin{eqnarray*}
M_j:=\{m\in \N\, |\, \psi(m)\geq j\}.
\end{eqnarray*}
The tuple $(M_1,...,M_{l_\psi})$ is called 
{\it standard covering of $\psi$}. It satisfies
(and is determined by this)
\begin{eqnarray*}
M_1\supset M_2\supset ...\supset M_{l_\psi}\neq\emptyset,
\quad \psi(m)=|\{j\in\{1,...,l_\psi\}\,|\, m\in M_j\}|.
\end{eqnarray*}

\bigskip
The most important result in this paper is Theorem \ref{t1.4}.

\begin{theorem}\label{t1.4}
(a) In the case of an isolated quasihomogeneous singularity,
each set $M_j$ satisfies condition (I) in Definition 
\ref{t2.6} (c).

\medskip
(b) 
$\Aut_{S^1}(H_{M_j},h_{M_j})=\{\pm h_{M_j}^k\,|\,k\in\Z\}.$
\end{theorem}

{\bf Proof:}
Part (b) follows from part (a) and Theorem \ref{t1.3}.
Part (a) follows from Lemma \ref{t2.7} and Theorem \ref{t6.6}.
Lemma \ref{t2.7} says that a finite set $M\subset \N$ 
which is compatible
with a tuple of excellent orders (Definition \ref{t2.4} (d))
satisfies condition (I). Theorem \ref{t6.6} says 
that in the case
of an isolated quasihomogeneous hypersurface singularity
each set $M_j$ is compatible with a certain tuple of 
excellent orders. It was proved in 
\cite[Theorem 12.1]{HM20-2}. It is cited in section 
\ref{s6}. \hfill$\Box$ 

\bigskip
{\it Problem 6} in \cite{HZ19} asked precisely for part
(a) of Theorem \ref{t1.4}. 
The notion of an excellent order and the notion of  
compatibility of a finite set $M\subset \N$ with a tuple
of excellent orders were first used in 
\cite[section 9]{HM20-2}. We found them when we searched
for a partial proof of Orlik's conjecture.
They are recalled in section \ref{s2}.

Theorem \ref{t1.4} is independent of Orlik's conjecture.
But it is useful only in the cases where Orlik's conjecture
holds. 

\begin{conjecture}\label{t1.5}
(Orlik's conjecture, \cite[conjecture 3.1]{Or72}) 
For any isolated quasihomogeneous singularity, 
there is an isomorphism
\begin{eqnarray*}
(H_{Mil},h_{Mil})\cong \bigoplus_{j=1}^{l_\psi}
(H_{M_j},h_{M_j}).
\end{eqnarray*}
\end{conjecture}

Our partial positive results \cite[Theorem 1.3]{HM20-2}
on this conjecture are recalled below in Theorem \ref{t3.4}
and right after it. 

In those isolated quasihomogeneous singularities, where 
Orlik's conjecture holds, $(H_{Mil},h_{Mil})$ decomposes (not 
uniquely in general) into Orlik blocks $(H_{M_j},h_{M_j})$,
and Theorem \ref{t1.4} says that the group of automorphisms
with eigenvalues in $S^1$ of each of these blocks
is the smallest possible, it is 
$\{\pm h_{M_j}^k\,|\, k\in\Z\}$.
This is very helpful also in determining the group of 
all automorphisms of the triple
$(H_{Mil},h_{Mil},L_{Mil})$ where $L_{Mil}:H_{Mil}\times
H_{Mil}\to\Z$ is the Seifert form. 
This group is important for period maps and Torelli 
conjectures for singularities (e.g. \cite{He11}). 

Problem 7 in \cite{HZ19} is a strengthening of problem 6
and is also solved by Lemma \ref{t2.7} and the results
in \cite[section 9]{HM20-2}. This is explained in section
\ref{s3}. The other problems 1--5 in \cite{HZ19} are
 discussed in the sections \ref{s5} and \ref{s6}.

\begin{notations}\label{t1.6}
$\N=\{1,2,3,..\}$, $\N_0=\{0,1,2,..\}$. 
Always $m\in\N$ and $n\in\N$. 
For $I\subset \R$, $\Z_I:=\Z\cap I$. 
The support of a map $\psi:\N\to\Q$ is the set
$\supp(\psi):=\{m\in\N\,|\, \psi(m)\neq 0\}$. 
We will consider only maps $\psi:\N\to\Q$ with
finite support. Most often, $\psi$ will be a map
$\psi:\N\to\N_0$ with finite support. 

The set of all prime numbers is $\PP\subset\N$. 
For $p\in\PP$ and $m\in\N$, $v_p(m)\in\N_0$ is the unique
number with $m=\prod_{q\in\PP}q^{v_q(m)}$. Also the 
projection 
$$\pi_p:\N\to \{n\in\N\,|\, p\not| n\},\quad m\mapsto
m\cdot p^{-v_p(m)},$$
will be used. 
\end{notations}

\section{Compatibility of a finite set of natural
numbers with a tuple of excellent orders}\label{s2}
\setcounter{equation}{0}

\noindent
This section proposes and discusses a 
condition for a finite set $M\subset \N$ of natural numbers
and a condition for a map $\psi:\N\to\N_0$
with finite support. 
They are given in Definition \ref{t2.4} (d) and (e).
They have a number of good properties, which will be 
given below in Lemma \ref{t2.5}, Lemma \ref{t2.7}, 
and the Theorems \ref{t3.2}, \ref{t3.4} and \ref{t3.6}. 
This is prepared by several definitions.
The Definitions \ref{t2.1}, \ref{t2.3} and \ref{t2.4},
Lemma \ref{t2.5}, and the Theorems \ref{t3.2}, \ref{t3.4}
and \ref{t3.6} are recalled from \cite[section 9]{HM20-2}.
Definition \ref{t2.6} is recalled from \cite{He20}.
Lemma \ref{t2.7} is new. It connects the condition in 
Definition \ref{t2.4} (d) for a finite set $M\subset\N$ 
with condition (I) (in Theorem \ref{t1.3} and Definition
\ref{t2.6} (c)).

\begin{definition}\label{t2.1} \cite[9.1]{HM20-2}
(a) An {\it excellent order} $\succ$ on a set
$\Z_{[0,s(\succ)]}$ for some bound $s(\succ)\in\N_0$
is a strict order (so transitive and for all 
$a,b\in \Z_{[0,s(\succ)]}$ either $a=b$ or
$a\succ b$ or $b\succ a$) which is determined by the
set 
\begin{eqnarray}\label{2.1}
S(\succ)&:=& \{k\in\Z_{[0,s(\succ)]}\, |\, k\succ 0\}
\end{eqnarray}
in the following way:
\begin{eqnarray}\label{2.2}
\left. \begin{array}{l}
\succ\textup{ equals }>\textup{ on }S(\succ)\cup\{0\},\\
\succ\textup{ equals }<\textup{ on }
\Z_{[0,s(\succ)]}-S(\succ).\end{array}\right\}
\end{eqnarray}
($S(\succ)=\emptyset$ is allowed.) 
The maximal element of $\Z_{[0,s(\succ)]}$ with  respect
to $\succ$ is called $s^+(\succ)$, so $s^+(\succ)\succ k$
for any other element $k\in\Z_{[0,s(\succ)]}$.

(b) The trivial excellent order is $\succ_0$
with $s(\succ_0):=0$, so it is the empty order on 
$\Z_{[0,s(\succ_0)]}=\{0\}$ 
(and, of course $S(\succ_0)=\emptyset$).

\medskip
(c) The tensor product 
of two excellent orders $\succ_1$ and $\succ_2$
is the excellent order $\succ_1\otimes\succ_2$ with 
\begin{eqnarray}\label{2.3}
s(\succ_1\otimes \succ_2)&:=&\max(s(\succ_1),s(\succ_2))
\quad\textup{and}\\
S(\succ_1\otimes \succ_2)&:=&
(S(\succ_1)\cup S(\succ_2))-(S(\succ_1)\cap S(\succ_2)).
\label{2.4}
\end{eqnarray}
\end{definition}

\begin{examples}\label{t2.2} \cite[9.2]{HM20-2}
(i) The excellent order $\succ_1$ with
$s(\succ_1)=7$ and $S(\succ_1)=\{6,4,1\}$ is given by
\begin{eqnarray*}
s^+(\succ_1)=6\succ_1 4\succ_1 1\succ_1 0 \succ_1 2
\succ_1 3\succ_1  5\succ_1 7.
\end{eqnarray*}
(ii)  The excellent order $\succ_2$ with
$s(\succ_2)=6$ and $S(\succ_2)=\{6,5,2,1\}$ is given by
\begin{eqnarray*}
s^+(\succ_2)=6\succ_2 5\succ_2 2\succ_2 1 
\succ_2 0\succ_2 3\succ_2 4.
\end{eqnarray*}
(iii) The excellent order $\succ_3:=(\succ_1\otimes \succ_2)$
for $\succ_1$ and $\succ_2$ in (i) and (ii) satisfies
$s(\succ_3)=7$, $S(\succ_3)=\{5,4,2\}$ and is given by
\begin{eqnarray*}
5\succ_3 4\succ_3 2\succ_3 0 \succ_3 1\succ_3 3\succ_3 6
\succ 7.
\end{eqnarray*}
(iv) For any excellent order $\succ$, the tensor product
with the trivial excellent order is $\succ$ itself, 
$\succ\otimes \succ_0=\succ$.
\end{examples}

\begin{definition}\label{t2.3} \cite[9.3]{HM20-2}
(a) A {\it path} in a finite directed graph $(V,E)$ 
(so $V$ is a finite non-empty set and $E\subset V\times V$)
is a  tuple $(v_1,...,v_l)$ for some $l\in\Z_{[2,\infty)}$ 
with $v_j\in V$ and $(v_j,v_{j+1})\in E$ for 
$j\in\Z_{[1,l-1]}$. It is a path from $v_1$ to $v_l$,
so with source $v_1$ and target $v_l$.

\medskip
(b) A finite directed graph $(V,E)$ {\it has a center}
$v_V\in V$ if it has no path from
any vertex to itself and if it has at least one path from
$v_V$ to any other vertex $v\in V$.
(The center is unique, which justifies the notation
$v_V$.)

\medskip
(c) Consider a tuple $(\succ_p)_{p\in P}$ of excellent 
orders for a finite set $P\subset \PP$ of prime numbers.
It defines a finite directed graph $(V,E_V)$ with center $v_V$
as follows. Its set $V=V((\succ_p)_{p\in P})$ of vertices
is the {\it quadrant} in $\N$ 
\begin{eqnarray}\label{2.5}
V&:=& \{\prod_{p\in P}p^{k_p}\,|\, 
k_p\in\Z_{[0,s(\succ_p)]}\}.
\end{eqnarray}
Its set of edges $E_V=E((\succ_p)_{p\in P})$ is the set
\begin{eqnarray}
E_V&:=& \bigcup_{p\in P}E_{V,p}\textup{ with }\label{2.6}\\
E_{V,p}&:=& \{(m_a,m_b)\in V\times V\,|\, 
\pi_p(m_a)=\pi_p(m_b), v_p(m_a)\succ_p v_p(m_b)\}.
\nonumber
\end{eqnarray}
The edges in $E_{V,p}$ are called $p$-edges.
So, the underlying undirected graph coincides with the
undirected graph which underlies the directed graph
$(V,E(V))$ in Definition \ref{t1.2}.
But the directions of edges may have changed.
The graph $(V,E_V)$ is obviously centered with center
\begin{eqnarray}\label{2.7}
v_V&=&\prod_{p\in P}p^{s^+(\succ_p)}.
\end{eqnarray}
\end{definition}

\begin{definition}\label{t2.4} \cite[9.4]{HM20-2}
(a) Let $\succ$ be an excellent order on the set
$\Z_{[0,s(\succ)]}$. A set $K\subset \N_0$ is 
{\it subset compatible} with $\succ$ if a bound
$k_K\in\Z_{[0,s(\succ)]}$ with 
\begin{eqnarray}\label{2.8}
K&=&\{k\in \Z_{[0,s(\succ)]}\, |\, k\succ k_K\}
\end{eqnarray}
exists or if $K=\Z_{[0,s(\succ)]}$. 
(The bound $k_K=s^+(\succ)$ gives $K=\emptyset$, 
which is allowed.)

\medskip
(b) For a finite non-empty set $M\subset\N$, let
\begin{eqnarray}\label{2.9}
\PP(M)&:=& \{p\in\PP\,|\, M\neq \pi_p(M)\}\\
&=& \{p\in \PP\,|\, \exists\ m\in M\textup{ with }v_p(m)>0\}
\nonumber
\end{eqnarray}
be the set of prime numbers which turn up as factors of some
numbers in $M$.

\medskip
(c) For a finite non-empty set $M\subset \N$, 
a prime number $p\in\PP$, and a number $m_0\in\pi_p(M)$, 
define the finite set $K_{M,p,m_0}\subset\N_0$ by 
\begin{eqnarray}\label{2.10}
\pi_p^{-1}(m_0)=\{m_0\cdot p^k\,|\, k\in K_{M,p,m_0}\}.
\end{eqnarray}

\medskip
(d) A finite non-empty set $M\subset\N$ is
{\it compatible} with a tuple $(\succ_p)_{p\in P}$
of excellent orders for a finite set $P\supset \PP(M)$ 
of prime numbers if 
\begin{eqnarray}\label{2.11}
M&\subset& V((\succ_p)_{p\in P})
\end{eqnarray}
and if for any prime number $p\in\PP(M)$
and any $m_0\in\pi_p(M)$ the set $K_{M,p,m_0}$
is subset compatible with $\succ_p$.
(So, here the excellent orders $\succ_p$ for 
$p\in \PP-\PP(M)$ are irrelevant.  But considering 
$P\supset \PP(M)$ instead of $P=\PP(M)$ will be useful.)

\medskip
(e) A map $\psi:\N\to\N_0$ with finite support 
$\supp(\psi)\subset\N$ is {\it compatible} with a tuple 
$(\succ_p)_{p\in P}$ of excellent orders for a finite set 
$P\supset \PP(\supp(\psi))$ of prime numbers if 
\begin{eqnarray}\label{2.12}
\supp(\psi)&\subset& V((\succ_p)_{p\in P})
\end{eqnarray}
and if for any edge $(m_a,m_b)\in E_V$
\begin{eqnarray}\label{2.13}
\psi(m_a)&\geq & \psi(m_b).
\end{eqnarray}

\medskip
(f) A {\it covering} of a map $\psi:\N\to \N_0$
with finite support is a tuple 
$(M_1,...,M_l)$ ($l\in\N_0$) 
of finite non-empty sets $M_j\subset\N$ with
\begin{eqnarray}\label{2.14}
\psi(m)&=& |\{j\in\{1,...,l\}\,|\, m\in M_j\}|
\textup{ for any }m\in\N.
\end{eqnarray} 
Here obviously $l\geq\max(\psi(m)\,|\, m\in\N_0)=:l_\psi$.
In the case $\supp(\psi)=\emptyset$ we have 
$l=0$ and an empty tuple.
The {\it standard covering} of $\psi$ is the tuple
$(M_1^{(st)},...,M_{l_\psi}^{(st)})$ with 
\begin{eqnarray}\label{2.15}
M_j^{(st)}=\{m\in \supp(\psi)\,|\, \psi(m)\geq j\}
\textup{ for }j\in\{1,...,l_\psi\}.
\end{eqnarray}
It is the unique covering with $M_1\supset ...\supset M_l$,
and it satisfies $M_1^{(st)}=\supp(\psi)$. 

\medskip
(g) Let $\psi:\N\to \N_0$ have finite support, 
let $P\supset \PP(\supp(\psi))$ be a finite set of prime
numbers, and let
$(\succ_p)_{p\in P}$ be a tuple of excellent orders
with \eqref{2.12}.
A covering $(M_1,...,M_l)$ of $\psi$ is called
{\it compatible} with $(\succ_p)_{p\in P}$
if each set $M_j$ is compatible with  $(\succ_p)_{p\in P}$. 
\end{definition}

The following lemma expresses the compatibility conditions
in Definition \ref{t2.4} (d) and (e) in a different way,
and it shows their relationship.
The proof is not difficult. It is given in \cite{HM20-2}.

\begin{lemma}\label{t2.5} \cite[9.5]{HM20-2}
(a) Let $M\subset\N$ be a finite non-empty set,
let $P\supset \PP(M)$ be a finite set of prime
numbers, 
and let $(\succ)_{p\in P}$ be a tuple of excellent orders
with \eqref{2.11}. (Recall the definition of $(V,E_V,v_V)$ 
in Definition \ref{t2.3} (c).) 
The following three conditions are equivalent:
\begin{list}{}{}
\item[(i)]
$M$ is compatible with $(\succ)_{p\in P}$.
\item[(ii)]
$(M,E_V\cap M\times M)$ is a directed graph with center $v_V$
(so $v_V\in M)$, and if $M$ contains the target of a path
in $(V,E_V)$, it contains all vertices in this path.
\item[(iii)]
If $m_b\in M$ and $(m_a,m_b)\in E_V$ then 
$m_a\in M$. 
\end{list}

\medskip
(b) Let $\psi:\N\to\N_0$ be a map with finite support,
let $P\supset \PP(\supp(\psi))$ be a finite set of prime
numbers, and let $(\succ_p)_{p\in P}$ be a tuple of 
excellent orders with \eqref{2.12}. 
The following three conditions are equivalent.
\begin{list}{}{}
\item[(i)] 
$\psi$ is compatible with $(\succ_p)_{p\in P}$.
\item[(ii)] 
$\psi$ has a covering $(M_1,...,M_l)$ which is compatible
with $(\succ_p)_{p\in P}$.
\item[(iii)]
The standard covering of $\psi$ is compatible with
$(\succ_p)_{p\in P}$.
\end{list}
\end{lemma}

Lemma \ref{t2.5} (a) allows to prove in Lemma \ref{t2.7} that
any finite non-empty set $M\subset\N$ which is
compatible with a tuple $(\succ_p)_{p\in P}$
of excellent orders satisfies condition (I)
in Definition \ref{t2.6} (c) and in Theorem \ref{t1.3}.
Before, we recall some necessary definitions from 
\cite{He20}.

\begin{definition}
\cite[Definition 1.1 and Theorem 1.2]{He20}\label{t2.6}
Let $M\subset \N$ be a finite non-empty set. 
Recall the definition of the graph $(M,E(M))$ in
Definition \ref{t1.2}.

\medskip
(a) For any prime number $p$, the components of the
graph $(M,E(M)-E_p(M))$ are called the 
{\it $p$-planes} of the graph. 
A $p$-plane is called a {\it highest} $p$-plane
if no $p$-edge ends at a vertex of the $p$-plane.
A $p$-edge from $m_1\in M$ to $m_2\in M$
is called a {\it highest} $p$-edge if no $p$-edge
ends at $m_1$.

\medskip
(b) Two properties $(T_p)$ and $(S_p)$ are defined for
any prime number $p$:
\begin{eqnarray}\label{2.16}
(T_p)&:& \textup{The graph }(M,E(M))\textup{ has only one
highest }p\textup{-plane}.\hspace*{1cm} \\
(S_p)&:& \textup{The graph }(M,E(M)-\{\textup{highest }p\textup{-edges}\}) \nonumber\\
&& \textup{ has only 1 or 2 components}.\label{2.17}
\end{eqnarray}

\medskip
(c) Condition (I) says: The graph $(M,E(M))$ is 
connected. It satisfies $(S_2)$. It satisfies $(T_p)$ 
for any prime number $p\geq 3$. 

\medskip
(d) Condition (II) says: The graph $(M,E(M))$ has two 
components $M_1$ and $M_2$. They are 2-planes.
They satisfy $(T_p)$ for any prime number $p\geq 3$. 
Furthermore,
\begin{eqnarray}\label{2.18}
\gcd(\lcm(M_1),\lcm(M_2))&\in &\{1,2\},\\
v_2(\lcm(M_2))&>&v_2(\lcm(M_1))\in\{0;1\}.\label{2.19}
\end{eqnarray}
\end{definition}

Lemma \ref{t2.7} is new. Part (b) is used in the proof of
Theorem \ref{t1.4}. 

\begin{lemma}\label{t2.7}
Let $M\subset\N$ be a finite non-empty set.

\medskip
(a) $(S_p)\Longrightarrow (T_p)$ for any prime number $p$
if the graph $(M,E(M))$ is connected.

\medskip
(b) Let $(\succ_p)_{p\in P}$ be a tuple of excellent orders
with which $M$ is compatible. 
Then the graph $(M,E(M))$ is connected.
$M$ satisfies $(S_p)$ for any prime number $p$.
$M$ satisfies condition (I).
\end{lemma}

{\bf Proof:}
(a) Suppose that the graph $(M,E(M))$ is connected
and satisfies $(S_p)$, but has at least two highest
$p$-planes.
Any highest $p$-plane is a component of the graph
$(M,E(M)-\{\textup{highest }p\textup{-edges}\})$.
Therefore this graph consists of exactly two
highest $p$-planes. But then the graph $(M,E(M))$ has
no $p$-edges and has the same two components,
so it is not connected, a contradiction.

(b) By Lemma \ref{t2.5} (a), $(M,E_V\cap M\times M)$ 
is a directed subgraph of $(V,E_V)$ with center $v_V$. 
Therefore it is connected. As the underlying undirected 
graphs of $(V,E(V))$ and $(V,E_V)$ coincide, also the graph 
$(M,E(M))$ is connected. 

Fix a prime number $p$. We will prove $(S_p)$ in the next
three paragraphs. Let $m_1\in M$ be arbitrary.
There is a path in $(V,E_V)$ from 
$\pi_p(v_V)\cdot p^{v_p(m_1)}$ to $m_1$ which consists
completely of edges which are not $p$-edges.
Because of Lemma \ref{t2.5} (a) (i)$\Rightarrow$(ii), it is 
a path in $(M,E_V\cap M\times M)$. 
Therefore the set $\{m\in M\,|\, v_p(m)=v_p(m_1)\}$
is a single $p$-plane in the graph $(M,E(M))$.
Furthermore, the set $\{\pi_p(v_V)\cdot p^k\,|\, 
k\in K_{M,p,\pi_p(v_V)}\}$ intersects each $p$-plane in
one element. 

Especially, there is only one highest $p$-plane,
and it intersects the set above in 
$\www v_V:=\pi_p(v_V)\cdot p^{k^{(max)}}$ where
$k^{(max)}:=\max(v_p(m)\,|\, m\in M)$.
This highest $p$-plane is one component of the graph 
$(M,E(M)-\{\textup{highest }p\textup{-edges}\})$.

The edges in $E(M)$ with both vertices in the set 
$\{\pi_p(v_V)\cdot p^k\,|\, 
k\in K_{M,p,\pi_p(v_V)}-\{k^{(max)}\}\}$
are all $p$-edges, but not highest $p$-edges. 
Therefore the second component (if $M$ does not consist only 
of one $p$-plane) is the union of all other $p$-planes.
The property $(S_p)$ holds. 

Because $(M,E(M))$ is connected and because of part (a),
$(M,E(M))$ satisfies condition (I).
\hfill$\Box$

\section{Tensor product and Thom-Sebastiani sum}\label{s3}
\setcounter{equation}{0}

\noindent
The following is a reformulation of {\it Problem 7}
in \cite{HZ19}:

{\it (a) Find a natural condition for a map $\psi:\N\to\N_0$
with finite support which implies that each set
$M_j^{st}\subset\N$ in the standard covering of $\psi$ 
satisfies condition (I) and which is stable under 
tensor product.

(b) Prove that in the case of an isolated quasihomogeneous
singularity the map $\psi:\N\to\N_0$ satisfies this 
condition, where 
$p_{ch,h_{Mil}}=\prod_{m\in\N}\Phi_m^{\psi(m)}$
is the characteristic polynomial of the monodromy.}

\bigskip
Problem 7 is motivated by the facts (which are discussed
in \cite{HZ19}) that condition (I) is not stable under
tensor product and that it is too difficult to prove
Theorem \ref{t1.4} directly. 

Here we propose as solution the natural condition 
that a tuple of excellent orders exists, 
with which $\psi$ is compatible (Definition 
\ref{t2.4} (e)).
The first half of part (a) is true because of Lemma 
\ref{t2.7} (b) and Lemma \ref{t2.5} (b). 
Part (b) is true because of Theorem \ref{t6.6}, which is 
\cite[Theorem 12.1]{HM20-2}. 
The second half of part (a) is also true. It follows 
from Theorem \ref{t3.2} below,
which is a main result of \cite{HM20-2}. 
Before, Definition \ref{t3.1}
explains what is meant here by {\it tensor product}.

\begin{definition}\label{t3.1}
(a) Let $f(t)=\prod_{i=1}^n(t-\kappa_i)\in\C[t]$ and 
$g(t)=\prod_{j=1}^n(t-\lambda_j)\in\C[t]$ be unitary
polynomials of degrees $n,m\geq 1$. Their tensor product
is the unitary polynomial
\begin{eqnarray}\label{3.1}
(f\otimes g)(t)&:=&\prod_{i=1}^n\prod_{j=1}^m
(t-\kappa_i\lambda_j)\in\C[t]
\end{eqnarray}
of degree $nm$.

\medskip
(b) Let $\psi_1:\N\to\N_0$ and $\psi_2:\N\to\N_0$
be two maps with finite supports $M=\supp(\psi_1)$
and $N=\supp(\psi_2)$. Their tensor product is the 
map $\psi_1\otimes\psi_2:\N\to\N_0$ with
\begin{eqnarray}\label{3.2}
\Bigl(\prod_{m\in M}\Phi_m^{\psi_1(m)}\Bigr)
\otimes \Bigl(\prod_{n\in N}\Phi_n^{\psi_2(n)}\Bigr)
= \prod_{l\in \N}\Phi_l^{(\psi_1\otimes\psi_2)(l)}.
\end{eqnarray}
It has also finite support. 
\end{definition}

The beginning of section 7 in \cite{HM20-2} 
(namely Definition 7.1, Lemma 7.2 and Definition 7.3)
gives more information on $\psi_1\otimes\psi_2$. 
The next Theorem \ref{t3.2} is a main result of 
\cite{HM20-2}. The long and difficult proof is built up from 
several intermediate results in the sections 7--9 
in \cite{HM20-2}. Its part (a) solves the second half
of part (a) of Problem 7 above.

\begin{theorem}\cite[Theorem 9.10]{HM20-2}\label{t3.2}
Let $\psi_1:\N\to\N_0$ and $\psi_2:\N\to\N_0$
be two maps with finite supports $M=\supp(\psi_1)$
and $N=\supp(\psi_2)$. Write $\psi_3:=\psi_1\otimes\psi_2$.
Denote $L:=\supp(\psi_3)$. 
It satisfies $\PP(L)\subset\PP(M)\cup\PP(N)$. 

Let $P\supset\PP(M)\cup\PP(N)$ be a finite set of
prime numbers. Let $(\succ_p^M)_{p\in P}$
and $(\succ_p^N)_{p\in P}$ be two tuples of excellent
orders such that $\psi_1$ is compatible with
$(\succ_p^M)_{p\in P}$ and $\psi_2$ is compatible
with $(\succ_p^N)_{p\in P}$. 
Write $\succ_p^L:=(\succ_p^M\otimes \succ_p^N)$
for any $p\in P$. 

\medskip
(a) $\psi_3$ is compatible with the tuple
$(\succ_p^L)_{p\in P}$ of excellent orders.

\medskip
(b) Let $(M_1^{(st)},...,M_{l_{\psi_1}}^{(st)})$,
$(N_1^{(st)},...,N_{l_{\psi_2}}^{(st)})$ and
$(L_1^{(st)},...,L_{l_{\psi_3}}^{(st)})$ be
the standard coverings of $\psi_1$, $\psi_2$
and $\psi_3$. Then
\begin{eqnarray}\label{3.3}
\Bigl(\bigoplus_{i=1}^{l_{\psi_1}}\Or(M_i^{(st)})\Bigr)
\otimes
\Bigl(\bigoplus_{j=1}^{l_{\psi_2}}\Or(N_j^{(st)})\Bigr)
\cong \bigoplus_{k=1}^{l_{\psi_3}}\Or(L_k^{(st)}),
\end{eqnarray}
so the tensor product of the sums of Orlik blocks on the left
hand side admits a standard decomposition into Orlik blocks.
\end{theorem}

Thom and Sebastiani proved a result which specializes
to the case of isolated quasihomogeneous singularities
as follows.

\begin{theorem}\cite{ST71}\label{t3.3}
Let $f=f(x_1,...,x_{n_f})$ and 
$g=g(x_{n_f+1},...,x_{n_f+n_g})$ be two isolated 
quasihomogeneous singularities in different variables.
Then their Thom-Sebastiani sum $f+g$ is an isolated 
quasihomogeneous singularity, and there is a canonical 
isomorphism
\begin{eqnarray}\label{3.4}
(H_{Mil},h_{Mil})(f+g)\cong 
(H_{Mil},h_{Mil})(f)\otimes (H_{Mil},h_{Mil})(g).
\end{eqnarray}
\end{theorem}

This theorem, Theorem \ref{t6.6} and part (b) of Theorem 
\ref{t3.2} imply the following main result of \cite{HM20-2}.

\begin{theorem}\cite[Theorem 1.3 (c)]{HM20-2}\label{t3.4}
Let $f=f(x_1,...,x_{n_f})$ and 
$g=g(x_{n_f+1},...,x_{n_f+n_g})$ be two isolated 
quasihomogeneous singularities in different variables,
which both satisfy Orlik's conjecture. 
Then also $f+g$ satisfies Orlik's conjecture.
\end{theorem}

Together with the proofs of Orlik's conjecture
for the chain type singularities 
\cite[Theorem 1.3 (a)]{HM20-2}
and for the cycle type singularities 
\cite[Theorem 1.3]{HM20-1},
it shows Orlik's conjecture for a large family of
isolated quasihomogeneous singularities, namely
for all invertible polynomials.

\begin{remark}\label{t3.5}
The search for Theorem \ref{t3.2} led us in the following way
to the excellent orders and the compatibility of a finite set 
$M\subset\N$  with a tuple of excellent orders. 
We found a condition {\it sdiOb-sufficient} 
\cite[Definition 7.3 (d)]{HM20-2} 
(short for {\it standard decomposition
into Orlik blocks}) for a pair $(M,N)$ of finite subsets
$M,N\subset\N$ which implies (and probably is equivalent to)
the property that the tensor product 
$(H_M,h_M)\otimes(H_N,h_N)$ of two Orlik blocks 
allows a standard decomposition into Orlik blocks 
\cite[Theorem 7.4]{HM20-2}.
Lemma 9.12 in \cite{HM20-2} says especially that the
following conditions (i) and (ii) are equivalent.
Here $M\subset\N$ is a finite non-empty set.
\begin{list}{}{}
\item[(i)]
For any $n_N\in\N$, the pair
$(M,\{n\in\N\,|\, n|n_N\})$ is sdiOb-sufficient.
\item[(ii)]
A tuple $(\succ_p)_{p\in \PP(M)}$ of excellent orders
exists such that $M$ is compatible with it.
\end{list}
\end{remark}

Another interesting result from \cite{HM20-2} which 
complements Theorem \ref{t3.2} (and which, in fact, is used 
in its proof) is the following.

\begin{theorem}\cite[Theorem 9.6]{HM20-2}\label{t3.6}
Let $(\succ_p)_{p\in P}$ be at tuple of excellent orders
for a finite set $P$ of prime numbers, and let 
$\psi:\N\to \N_0$ be a map with finite support which
is compatible with $(\succ_p)_{p\in P}$. 
Let $(M_1^{(1)},...,M_{l_1}^{(1)}))$ and $(M_1^{(2)},...,
M_{l_2}^{(2)})$ be two coverings of $\psi$ which
are both compatible with $(\succ_p)_{p\in P}$.
Then the corresponding sums of Orlik blocks
are isomorphic,
\begin{eqnarray}\label{3.5}
\bigoplus_{i=1}^{l_1}\Or(M_i^{(1)})&\cong &
\bigoplus_{j=1}^{l_2}\Or(M_j^{(2)}),
\end{eqnarray}
and $l_1=l_2=l_\psi(:=\max(\psi(m)\,|\, m\in\N))$.
\end{theorem}

\section{Maps on the group ring generated by unit roots}
\label{s4}
\setcounter{equation}{0}

\noindent
This section collects notations and classical formulas
which will be needed in the later sections.

\begin{notations}\label{t4.1}
(a) 
Denote by $\mu(\C)\subset S^1$ the group of all unit roots.
Denote by $\Q[\mu(\C)]$ and $\Z[\mu(\C)]$ the group rings
of elements $\sum_{j=1}^lb_j[\zeta_j]$ where 
$\zeta_j\in\mu(\C)$ and where $b_j\in\Q$ respectively $\Z$,
with multiplication $[\zeta_1][\zeta_2]=[\zeta_1\zeta_2]$.
The unit element is $[1]$. The trace and the degree of
an element $b=\sum_{j=1}^lb_j[\zeta_j]\in\Q[\mu(\C)]$ are
\begin{eqnarray*}
\tr(b):=\sum_{j=1}^lb_j\zeta_j\in\C,\quad
\deg(b):=\sum_{j=1}^lb_j\in\Q.
\end{eqnarray*}
The trace map $\tr:\Q[\mu(\C)]\to\C$ and the degree map
$\deg:\Q[\mu(\C)]\to\Q$ are ring homomorphisms.
For $k\in\N_0$, the $k$-th Lefschetz number $L_k(b)$
of an element $b=\sum_{j=1}^lb_j[\zeta_j]$ is
\begin{eqnarray}\label{4.1}
L_k(b):=\sum_{j=1}b_j\zeta_j^k
=\tr\left(\sum_{j=1}b_j[\zeta_j^k]\right)\in\C.
\end{eqnarray}

\medskip
(b) The divisor of a unitary polynomial 
$f=\prod_{i=1}^n(t-\kappa_i)$ with $\kappa_i\in\mu(\C)$ is 
\begin{eqnarray*}
\divis f:=\sum_{i=1}^n[\kappa_i]\in\Z[\mu(\C)].
\end{eqnarray*}
The divisor of the constant polynomial 1 is $\divis 1=0$.

\medskip
(c) The order $\ord(\zeta)\in\N$ of a unit root 
$\zeta\in\mu(\C)$ is the minimal number $m\in\N$ with
$\zeta^m=1$. For $m\in\N$, the $m$-th cyclotomic polynomial
is 
\begin{eqnarray*}
\Phi_m:=\prod_{\zeta:\,\ord(\zeta)=m}(t-\zeta).
\end{eqnarray*}
Define 
\begin{eqnarray}\label{4.2}
\Lambda_m:=\divis(t^m-1),\quad \Psi_m:=\divis\Phi_m.
\end{eqnarray}
Then $\lambda_1=\Psi_1=[1]$. 
The M\"obius function $\mu_{Moeb}$ is \cite{Ai79}
\begin{eqnarray*}
\mu_{Moeb}:\N&\to&\{0,1,-1\},\\
m&\mapsto& \left\{\begin{array}{ll}
(-1)^r & \textup{if }m=p_1\cdot ...\cdot p_r\textup{ with }
p_1,...,p_r\\
 & \textup{different prime numbers},\\
0& \textup{else}  
\end{array}\right. \nonumber
\end{eqnarray*}
(here $r=0$ is allowed, so $\mu_{Moeb}(1)=1$).
\end{notations}

The next lemma collects well known facts.

\begin{lemma}\label{t4.2}
(a) Let $f=\prod_{i=1}^n(t-\kappa_i)$ and
$g=\prod_{j=1}^m(t-\lambda_j)$ be unitary polynomials with 
$\kappa_i$ and $\lambda_j\in\mu(\C)$. The tensor product 
$(f\otimes g)(t)=\prod_{i=1}^n\prod_{j=1}^m
(t-\kappa_i\lambda_j)$ in \eqref{3.1} satisfies
\begin{eqnarray*}
\divis (f\otimes g) =(\divis f)\cdot(\divis g).
\end{eqnarray*}

(b) The cyclotomic polynomial $\Phi_m$ 
is in $\Z[t]$, it has degree $\varphi(m)$, 
and it is irreducible in $\Z[t]$ and $\Q[t]$.
\begin{eqnarray}\label{4.3}
t^n-1&=&\prod_{m|n}\Phi_m(t),\hspace*{2cm} 
\Lambda_n=\sum_{m|n}\Psi_m,\\
\Phi_m&=&\prod_{m|n}(t^n-1)^{\mu_{Moeb}(\frac{m}{n})},\quad 
\Psi_m=\sum_{n|m}\mu_{Moeb}(\frac{m}{n})\cdot\Lambda_n.
\label{4.4}
\end{eqnarray}
The traces of $\Lambda_m$ and $\Psi_m$ are 
\begin{eqnarray}\label{4.5}
\tr\Lambda_m&=&\left\{\begin{array}{ll}
1&\textup{if }m=1\\
0&\textup{if }m\geq 2,\end{array}\right. \\
\tr\Psi_m&=&\mu_{Moeb}(m).\label{4.6}
\end{eqnarray}
The product $\Lambda_m\Lambda_n$ is
\begin{eqnarray}\label{4.7}
\Lambda_m \Lambda_n &=&
\gcd(m,n)\cdot \Lambda_{\lcm(m,n)}.
\end{eqnarray}
For the product $\Psi_m\Psi_n$ see 
\cite[(2.19)--(2.21)]{HZ19} or \cite[Lemma 7.2]{HM20-2}.

\medskip
(c) Let $\psi:\N\to\Q$ be a map with finite support,
and consider the element 
\begin{eqnarray}\label{4.8}
b:=\sum_{m\in\N}\psi(m)\cdot\Psi_m\in \Q[\mu(\C)].
\end{eqnarray}
There is a unique map $\chi:\N\to\Q$ with finite
support and
\begin{eqnarray}\label{4.9}
b=\sum_{n\in\N}\chi(n)\cdot\Lambda_n.
\end{eqnarray}
The four data $\psi$, $b$, $\chi$ and $(L_k(b))_{k\in\N_0}$
determine one another, by \eqref{4.8}--\eqref{4.9} and 
\begin{eqnarray}
\psi(m)=\sum_{n:\, m|n}\chi(n),\quad 
\chi(n)=\sum_{m:\, n|m}\psi(m)\cdot\mu_{Moeb}
(\frac{m}{n}),\label{4.10}\\
L_k(b)= \sum_{n:\, n|k}n\chi(n),\quad
n\chi(n)=\sum_{k|n}L_k(b)\cdot\mu_{Moeb}(\frac{n}{k}).
\label{4.11}
\end{eqnarray}
Denote $d_\chi:=\lcm(n\,|\, n\in\supp(\chi))$. 
The Lefschetz numbers $L_k(b)$ are rational.
They are determined
by the values $L_k(b)$ for $k\in\N$ with $k| d_\chi$, because of the (extended) periodicity property 
\begin{eqnarray}\label{4.12}
L_k(b)=L_{\gcd(k,d_\chi)}(b)\qquad\textup{for any }
k\in\N.
\end{eqnarray}
\end{lemma}

{\bf Proof:}
(a) Trivial. (b) The statements on $\Phi_m$ are classical.
The formulas \eqref{4.4} follow from the formulas \eqref{4.3}
by Moebius inversion \cite{Ai79}.
The trace of $\Lambda_m$ is obvious, the trace of $\Psi_m$
follows with \eqref{4.4}. The product $\Lambda_m\Lambda_n$
is obvious.

(c) Formula \eqref{4.10} is a consequence of \eqref{4.3} 
and \eqref{4.4}. The first formula in \eqref{4.11}
can be seen as follows, 
\begin{eqnarray*}
L_k(\Lambda_n)
&=&\tr\left(\sum_{a=0}^{n-1}[e^{2\pi i ak/n}]\right)\\
&=&\tr\left(\gcd(k,n)\Lambda_{n/\gcd(k,n)}\right)
=\left\{\begin{array}{ll}
n&\textup{ if }n|k,\\ 0&\textup{ if }n\not| k,\end{array}
\right.  \\
L_k(b)&=&\sum_{n\in\N}\chi(n)\cdot L_k(\Lambda_n)
=\sum_{n:\, n|k}n\chi(n).
\end{eqnarray*}
The first formula in \eqref{4.11} gives the second formula 
in \eqref{4.11} by Moebius inversion \cite{Ai79}. 
It also implies $L_k(b)\in\Q$ and the periodicity property 
\eqref{4.12}. \hfill$\Box$

\section{Basic properties of weight systems and induced
objects}\label{s5}
\setcounter{equation}{0}

\noindent
In this section, problem 1 in \cite{HZ19} will be solved
as part of Lemma \ref{t5.4}. It is elementary.
Before, notations from \cite[section 3]{HZ19} are recalled.
In this section, we fix a number $n\in\N$ and denote
$N:=\{1,2,...,n\}$.

\begin{definition}\label{t5.1}
(a) A weight system is a tuple 
$(v_1,...,v_n;d)\in(\Q_{>0})^{n+1}$ with $v_i<d$.
Another weight system is {\it equivalent}, if the second
one has the form $q\cdot(v_1,...,v_n;d)$ for some 
$q\in\Q_{>0}$. A weight system is {\it integer} if
$(v_1,...,v_n;d)\in\N^{n+1}$. It is {\it reduced} if
it is integer and $\gcd(v_1,...,v_n,d)=1$.
It is {\it normalized} if $d=1$. 

Any equivalence class contains a unique reduced weight system
and a unique normalized weight system.
From now on, the letters $(v_1,...,v_n;d)$ will be 
reserved for an integer weight system, and the letters
$(w_1,...,w_n;1)$ will be the equivalent normalized 
weight system, i.e. $w_i=v_i/d$. 

(b) Let $(v_1,...,v_n;d)$ be an integer weight system
(not necessarily reduced). For $J\subset N$ and $k\in\N_0$
define
\begin{eqnarray}\label{5.1}
(\Z^J)_k&:=&\{\alpha\in\Z^n\,|\, \alpha_i=0\textup{ for }
i\notin J,\ \sum_{i\in J}\alpha_i=k\},\\
(\N_0^J)_k&:=& \N_0^n\cap (\Z^J)_k.\nonumber
\end{eqnarray}

(c) Two conditions for an integer weight system
$(v_1,...,v_n;d)$ (or, equivalently, for the 
equivalent normalized weight system $(\frac{v_1}{d},...,
\frac{v_n}{d};1)$) are defined as follows:
\begin{eqnarray*}
\overline{(C2)}:&& 
\forall\ J\subset N\textup{ with }J\neq \emptyset
\quad \exists\ K\subset N \\
&& \textup{with }|K|=|J|\textup{ and }
\forall\ k\in K\ (\Z^J)_{d-v_k}\neq\emptyset.\\
(C2):&& 
\forall\ J\subset N\textup{ with }J\neq \emptyset
\quad \exists\ K\subset N \\
&& \textup{with }|K|=|J|\textup{ and }
\forall\ k\in K\ (\N_0^J)_{d-v_k}\neq\emptyset.
\end{eqnarray*}
\end{definition}

\begin{remarks}\label{t5.2}
(i) Of course $(C2)$ implies $\oooo{(C2)}$. 
Lemma 2.1 in \cite{HK12} gives four conditions
$(C1)$, $(C1)'$, $(C2)'$ and $(C3)$ which are equivalent
to $(C2)$. They all have versions 
$\oooo{(C1)}$, $\oooo{(C1)'}$, $\oooo{(C2)'}$ and 
$\oooo{(C3)}$, which are equivalent to $\oooo{(C2)}$. 
Theorem \ref{t6.1} below shows the relevance of $(C2)$.

\medskip
(ii) Let $(v_1,...,v_n;d)$ be an integer weight system. 
For $J\subset N$ with $J\neq\emptyset$ define the semigroup
\begin{eqnarray*}
SG(J)&:=&\sum_{j\in J}\N_0\cdot v_j\subset\N_0,
\end{eqnarray*}
and observe $\sum_{j\in J}\Z\cdot v_j=\Z\cdot
\gcd(v_j\,|\, j\in J)$. Therefore 
\begin{eqnarray}\label{5.2}
(\N_0^J)_k\neq\emptyset&\iff& k\in SG(J),\\
(\Z^J)_k\neq\emptyset&\iff& \gcd(v_j\,|\, j\in J)|k.
\label{5.3}
\end{eqnarray}
$\oooo{(C2)}$ is equivalent to the condition
\begin{eqnarray*}
(GCD)&& \forall\ J\subset N\ \gcd(v_j\,|\, j\in J)
\textup{ divides at least }|J|\\
&&\textup{of the numbers }d-v_k\textup{ for }k\in N.
\end{eqnarray*}
\end{remarks}

\begin{definition}\label{t5.3}
Let $({\bf v},d)=(v_1,...,v_n,d)\in\N^{n+1}$ be an 
integer weight system. 

(a) Define unique numbers $s_1,...,s_n,t_1,...,t_n\in\N$ by 
\begin{eqnarray}\label{5.4}
\frac{v_i}{d}=\frac{s_i}{t_i}\quad\textup{and}\quad
\gcd(s_i,t_i)=1.
\end{eqnarray}
They depend only on the normalized weight system
${\bf w}=(w_1,...w_n)=(\frac{v_1}{d},...,\frac{v_n}{d})$.
Define 
\begin{eqnarray}\label{5.5}
d_{\bf w}:=\lcm(t_j\, |\,  j\in N).
\end{eqnarray}
Of course $d_{\bf w}|d$. If $({\bf v},d)$ is reduced
and $\gcd(v_1,...,v_n)|d$ (which holds for example
if $\oooo{(C2)}$ holds), then 
$\gcd(v_1,...,v_n)=1$ and then $d_{\bf w}=d$.

(b) For $k\in\N$ define
\begin{eqnarray}\label{5.6}
M(k)&:=&\{j\in N\, |\, t_j|k\},\\
\textup{and }\mu(k)&:=&\prod_{j\in M(k)}(\frac{1}{w_j}-1)
=\prod_{j\in M(k)}\frac{d-v_j}{v_j}\in\Q_{>0}\label{5.7}
\end{eqnarray}
(the empty product is by definition $1$).

(c) Define a quotient of polynomials
\begin{eqnarray}\label{5.8}
\rho_{({\bf v},d)}(t)&:=&t^{v_1+...+v_n}\cdot 
\prod_{j=1}^n \frac{t^{d-v_j}-1}{t^{v_j}-1}
\in\Q(t)
\end{eqnarray}
and an element of $\Q[\mu(\C)]$
\begin{eqnarray}\label{5.9}
D_{\bf w}&:=& \prod_{j=1}^n 
\left(\frac{1}{s_j}\Lambda_{t_j}-\Lambda_1\right) \in 
\Q[\mu(\C)].
\end{eqnarray}
\end{definition}

\begin{lemma}(\cite[Lemma 3.7]{HZ19} for (a)--(c))\label{t5.4}
Let $({\bf v};d)=(v_1,...,v_n;d)\in\N^{n+1}$ be an integer
weight system.

(a) Then $M(k)=M(\gcd(k,d_{\bf w}))$ and 
$\mu(k)=\mu(\gcd(k,d_{\bf w}))$. 

(b) The Lefschetz numbers
$L_k(D_{\bf w})$ are 
\begin{eqnarray}\label{5.10}
L_k(D_{\bf w})&=& (-1)^{n-|M(k)|}\mu(k)\\
&=& L_{\gcd(k,d_{\bf w})}(D_{\bf w}).\label{5.11}
\end{eqnarray}

(c) $({\bf v};d)$ satisfies $\oooo{(C2)}\iff
\rho_{({\bf v};d)}\in\Z[t]$.

(d) Suppose that $({\bf v};d)$ satisfies $\oooo{(C2)}$.
Write
\begin{eqnarray*}
\rho_{({\bf v};d)}=\sum_{\alpha\in d^{-1}\Z}
\sigma(\alpha)\cdot t^{d\cdot\alpha}
\end{eqnarray*}
with a unique map $\sigma:d^{-1}\Z\to\Z$
with finite support. Then
\begin{eqnarray}\label{5.12}
D_{\bf w}&=& \sum_{\alpha\in d^{-1}\Z}
\sigma(\alpha)\cdot [e^{2\pi i \alpha}].
\end{eqnarray}
and $\mu(k)\in\N$. 
\end{lemma}

{\bf Proof:}
(a) All $t_j$ divide $d_{\bf w}$.
Therefore $t_j|k\iff t_j|\gcd(k,d_{\bf w})$. 

(b) For completeness sake, the calculation in \cite{HZ19} 
is copied.
\begin{eqnarray*}
L_k(D_{\bf w})&=& \tr\left(\prod_{j=1}^n 
\left(\frac{\gcd(k,t_j)}{s_j}
\Lambda_{t_j/\gcd(k,t_j)}-\Lambda_1\right)\right)\\
&=& \prod_{j=1}^n \left( \frac{\gcd(k,t_j)}{s_j}\cdot
\tr\left(\Lambda_{t_j/\gcd(k,t_j)}\right)-1\right) \\
&=& \prod_{j\in M(k)}\left(\frac{t_j}{s_j}-1\right) 
\cdot \prod_{j\notin M(k)}(-1)\\
&=& (-1)^{n-|M(k)|}\cdot \mu(k).
\end{eqnarray*}
Formula \eqref{5.11} follows from part (a).

(c) The conditions $\oooo{(C2)}$ and $(GCD)$ in Remark
\ref{t5.2} are equivalent. Condition $(GCD)$ says that any
cyclotomic polynomial in the denominator of 
$\rho_{({\bf v};d)}$ turns up with at least the same
multiplicity in the numerator. Therefore $(GCD)$ is 
equivalent to $\rho_{({\bf v};d)}\in\Z[t]$. 

(d) The formulas \eqref{5.10} and \eqref{5.12} imply
$\mu(k)\in\N$. For \eqref{5.12}, it is because of Lemma
\ref{t4.2} (c) sufficient to show
\begin{eqnarray}\label{5.13}
L_k(D_{\bf w})&=& L_k\left(\sum_{\alpha\in d^{-1}\Z}
\sigma(\alpha)\cdot [e^{2\pi i\alpha}]\right)
\quad\textup{for all }k\in\N_0.
\end{eqnarray}
The right hand side of \eqref{5.13} is
\begin{eqnarray*}
&&\tr\left(\sum_{\alpha\in d^{-1}\Z}
\sigma(\alpha)\cdot [e^{2\pi ik\alpha}]\right)
= \sum_{\alpha\in d^{-1}\Z}
\sigma(\alpha)\cdot e^{2\pi ik\alpha}
= \rho_{({\bf v};d)}(e^{2\pi i k/d})\\
&=& \lim_{t\mapsto e^{2\pi i k/d}}
\left(\prod_{j\in M(k)}\frac{t^d-t^{v_j}}{t^{v_j}-1}
\cdot \prod_{j\notin M(k)}\frac{t^d-t^{v_j}}{t^{v_j}-1}
\right)\\\
&=& \left(\prod_{j\in M(k)}\frac{d-v_j}{v_j}\right)
\cdot (-1)^{n-|M(k)|} =(-1)^{n-|M(k)|}\mu(k)
=L_k(D_{\bf w}).
\end{eqnarray*}
Here we used $j\notin M(k)\iff \frac{k}{d}v_j\notin\Z
\iff (e^{2\pi i k/d})^{v_j}\neq 1$. \hfill$\Box$

\begin{remark}\label{t5.5}
{\it Problem 1} in \cite{HZ19} asked whether 
\eqref{5.12} holds. 
In \cite{HZ19}, the first author had missed the calculation
in the proof of part (d).
\end{remark}

\section{Isolated quasihomogeneous singularities}\label{s6}
\setcounter{equation}{0}

\noindent
Finally we come to the discussion of weight systems which
allow isolated quasihomogeneous singularities and of
associated objects.
The following Theorems \ref{t6.1} and \ref{t6.4} are basic.
As in section \ref{s5}, $n\in\N$ and $N=\{1,...,n\}$
are fixed.

\begin{theorem} (Kouchnirenko \cite[Remarque 1.13 (ii)]{Ko76})
\label{t6.1}
Let $({\bf v};d)=(v_1,...,v_n;d)$ be an integer weight system.
The following conditions are equivalent.
\begin{list}{}{}
\item[(IS3):] \quad There exists a quasihomogeneous polynomial 
$f$ with weight system $({\bf v};d)$ and 
an isolated singularity at 0.
\item[(IS3)':] \quad A generic quasihomogeneous polynomial 
with weight system $({\bf v};d)$ has an isolated 
singularity at 0.
\item[$(C2)$:] \quad The weight system $({\bf v};d)$ 
satisfies condition $(C2)$.
\end{list}
\end{theorem}

\begin{remarks}\label{t6.2}
(i) The reference \cite[Remarque 1.13 (ii)]{Ko76} does
not give a detailed proof, but the reference \cite{Ko77}
(in russian) does. The theorem was rediscovered and reproved
several times, in \cite{OR76}, \cite{KS92} and \cite{Fl00}.
More general results are given in \cite{Sh79} (without proof)
and \cite{Wa96} (with proof). 
The paper \cite{HK12} discussed the history of Theorem 
\ref{t6.1}, but it missed the reference \cite{Fl00}. 

\medskip
(ii) Lemma 2.1 in \cite{HK12} shows combinatorially the
equivalence of (C2) with several other conditions
(C1), (C1)', (C2)' and (C3). Some of the references above
prove the equivalence of (IS3) with some other of these
conditions. 

\medskip
(iii) The implication (IS3)$\Rightarrow$(C1) was already 
shown in \cite{Sa71}. 
\end{remarks}

\begin{remark}\label{t6.3}
The paper \cite{HK12} made good use of the part of
condition (C2) which concerns subsets $J\subset N$
with $|J|=1$. {\it Problem 3} in \cite{HZ19} asked about 
making good use of the conditions in (C2) for sets
$J\subset N$ with $|J|\geq 2$. Of course,
the proofs of Theorem \ref{t6.1} in the references
in Remark \ref{t6.2} (i) make such good use, as they
need the full condition (C2). 
We do not have other solutions of problem 3 in \cite{HZ19}.
\end{remark}

\begin{theorem}\label{t6.4}
Let $({\bf v};d)=(v_1,...,v_n;d)$ be an integer 
weight system which satisfies condition (C2),
and let ${\bf w}=(\frac{v_1}{d},...,\frac{v_n}{d};1)$
be the equivalent normalized weight system. 

\medskip
(a) (Milnor and Orlik \cite{MO70})
$D_{\bf w}$ is the divisor of the characteristic polynomial
of the monodromy on the Milnor lattice of an isolated
quasihomogeneous singularity with weight system
$({\bf v};d)$. Especially $D_{\bf w}\in\N_0[\mu(\C)]$.

\medskip
(b) (Many people, e.g. \cite{AGV85})
$t^{-v_1-...-v_n}\rho_{({\bf v};d)}$ is the generating 
function of the weighted degrees of the Jacobi algebra
$$\C[x_1,...,x_n]/\Bigl(\frac{\paa f}{\paa x_1},...,
\frac{\paa f}{\paa x_n}\Bigr)$$
of an isolated quasihomogeneous singularity $f$ with
weight system $({\bf v};d)$. 
Especially $\rho_{({\bf v};d)}\in\N_0[t]$. 
\end{theorem}

\begin{remarks}\label{t6.5}
(i) By Lemma \ref{t5.4} (d), $\rho_{({\bf v};d)}$ and
$D_{\bf w}$ are related by
\begin{eqnarray}\label{6.1}
\rho_{({\bf v};d)}=\sum_{\alpha\in d^{-1}\Z}\sigma(\alpha)
\cdot t^{d\cdot \alpha},\qquad
D_{\bf w}=\sum_{\alpha\in d^{-1}\Z}\sigma(\alpha)\cdot
[e^{2\pi i \alpha}],
\end{eqnarray} 
where $\sigma:d^{-1}\Z\to\N_0$ has finite support. 

\medskip
(ii) If we write $\rho_{({\bf v};d)}=\sum_{i=1}^\mu
t^{d\cdot\alpha_i}$ for suitable $\alpha_i\in\Q$,
then these numbers $(\alpha_1,...,\alpha_\mu)$ are the
{\it exponents} of the singularity, and 
$e^{2\pi i\alpha_1},...,e^{2\pi i\alpha_\mu}$ are the
eigenvalues of the monodromy. 

\medskip
(iii) {\it Problem 2} in \cite{HZ19} asked for a 
combinatorial proof of the fact 
$\rho_{({\bf v};d)}\in\N_0[t]$ for a weight system
$({\bf v};d)$ with condition (C2). 
(The other two parts of Problem 2 in \cite{HZ19} would 
follow from this and from Lemma \ref{t5.4} (d).)
Theorem \ref{t6.4} (b), which was rediscovered by many people,
contains this statement. The classical proof uses that
the tuple of partial derivatives 
$(\frac{\paa f}{\paa x_1},...,\frac{\paa f}{\paa x_n})$
is a regular sequence. We do not know a proof which 
does not use the Jacobi algebra and this regular sequence.
So we do not have a solution to problem 2 in \cite{HZ19}. 
\end{remarks}

The following theorem was proved in \cite{HM20-2} as a step
in a partial proof of Orlik's conjecture \ref{t1.5}.
It is also a step in the proof of Theorem \ref{t1.4}. 
Its proof in \cite{HM20-2} uses the Theorems \ref{t6.1} 
and \ref{t6.4} (a) and especially that certain subsystems 
${\bf w}^{(p)}$ of ${\bf w}$ (for certain prime numbers
$p$) also satisfy condition (C2) and thus that 
$D_{{\bf w}^{(p)}}\in\N_0[\mu(\C)]$ holds.

\begin{theorem}\label{t6.6}
Consider a normalized weight system ${\bf w}=(w_1,...,w_n;1)$
which satisfies condition (C2), and consider the numbers 
$s_i,t_i\in\N$ with $w_i=\frac{s_i}{t_i}$ and 
$\gcd(s_i,t_i)=1$. Write $D_{\bf w}=\sum_{m\in M_{\bf w}}
\psi_{\bf w}(m)\cdot\Psi_m$ where the map 
$\psi_{\bf w}:M_{\bf w}\to\N_0$ has finite support
$M_{\bf w}\subset\N$. 
The map $\psi_{\bf w}$ is compatible with the tuple
$(\succ_p^{\bf w})_{p\in\PP(M_{\bf w})}$ of excellent orders
which is defined as follows,
\begin{eqnarray}\label{6.2}
s(\succ_p^{\bf w})&:=& \max(v_p(m)\,|\, m\in M_{\bf w}),\\
S(\succ_p^{\bf w})&:=& \{k\in\Z_{[1,s(\succ_p)]}\,|\, 
|\{j\in N\,|\, p^k|t_j\}|\textup{ is odd.}\}\label{6.3}
\end{eqnarray}
\end{theorem}

Finally, we solve {\it problem 4} in \cite{HZ19} in
Remark \ref{t6.8} and discuss {\it problem 5} in
\cite{HZ19} in Remark \ref{t6.9}.

\begin{remark}\label{t6.7}
An isolated quasihomogeneous singularity 
$f\in\C[x_1,...,x_n]$ has {\it multiplicity $\geq 3$}
if all monomials in it have degree $\geq 3$
(the standard degree, not the weighted degree). 
A result in \cite{Sa71} says first that the weight
system ${\bf w}$ of any isolated quasihomogeneous singularity 
with multiplicity $\geq 3$ is unique and satisfies 
$w_i\in(0,\frac{1}{2})$ for all $i$, 
and second that any quasihomogeneous
singularity can be transformed by a coordinate change
to a Thom-Sebastiani sum $f(x_1,...,x_k)+x_{k+1}^2+...+x_n^2$ 
of an isolated quasihomogeneous singularity $f$ of 
multiplicity $\geq 3$ and the $A_1$-singularity
$x_{k+1}^2+...+x_n^2$. 
Because of Theorem \ref{t3.3}, one can often restrict
to (isolated quasihomogeneous) singularities of multiplicity
$\geq 3$. 
\end{remark}

\begin{remark}\label{t6.8}
Let $f$ be an isolated quasihomogeneous singularity
of multiplicity $\geq 3$ with weight system
${\bf w}=(w_1,...,w_n;1)$. 
If $n\leq 3$, $(C2)\iff\oooo{(C2)}$.
This was first proved in \cite[Theorem 3]{Sa87}.
A short combinatorial proof is given in \cite[Lemma 2.5]
{HK12}. For $n\geq 4$, $(C2)$ is stronger than $\oooo{(C2)}$.
Though only one weight system ${\bf w}$
with $n=4$ which satisfies $\oooo{(C2)}$, but not $(C2)$, 
is published, the weight system $(58,33,24,1;265)$. 
It is due to Ivlev \cite{AGV85}. 
It is also discussed in \cite[Example 5.1]{HZ19}. 
There {\it problem 4} is posed. It asks for other 
weight systems ${\bf w}$ with $n=4$ which satisfy 
$\oooo{(C2)}$, but not $(C2)$.

We used the software PARI/GP \cite{PARI19} in order to 
check all reduced integer weight systems 
$({\bf v};d)=(v_1,v_2,v_3,v_4;d)$ with 
$\frac{d}{2}>v_1\geq v_2\geq v_3\geq v_4$ and $d\leq 360$
(reduced means $\gcd(v_1,v_2,v_3,v_4,d)=1$).
We found 654077 such weight systems which satisfy 
$\oooo{(C2)}$.
Only 103 of them do not satisfy $(C2)$. 
Up to $d=200$, there are 185013 weight systems as above 
with $\oooo{(C2)}$, and only 23 of them do not satisfy $(C2)$. 
Table \ref{table 1} lists these 23 weight systems 
$(v_1,v_2,v_3,v_4;d)$ in the lexicographic ordering 
with respect to $(d,v_1,v_2,v_3,v_4)$.
It also gives the Milnor number $\mu$, the number $L$ of the
weight system in the lexicographic ordering with respect to 
$(d,v_1,v_2,v_3,v_4)$ of all weight systems as above which
satisfy $\oooo{(C2)}$, and a tuple 
$[a_1,a_2,a_3,a_4,a_5,a_6]\in\{0,1,2,3,4\}^6$ which says
the following. Write 
$$(J_1,J_2,J_3,J_4,J_5,J_6)=(\{1,2\},\{1,3\},\{1,4\},
\{2,3\},\{2,4\},\{3,4\}).$$
Then $a_i=|\{k\in N\,|\, 
(\N_0^{J_i})_{d-v_k}\neq \emptyset\}|$. 
That $(C2)$ is not satisfied implies that for at least 
one $i\in\{1,2,3,4,5,6\}$ the number $a_i$ is $\leq 1$. 
\begin{table}
\[\begin{array}{r|r|r|r|c}
(v_1,v_2,v_3,v_4) & d  & \mu & L & [a_1,a_2,a_3,a_4,a_5,a_6]\\
\hline 
(27,16,10,1)  & 81  & 4615 & 22598 & [2, 2, 4, 1, 4, 4]\\
(29,14,12,3) & 87  & 1825 & 26464 & [2, 2, 4, 1, 4, 2]\\
(31,18,10,3) & 93  & 2075 & 31738 & [2, 2, 4, 1, 2, 4]\\
(35,26,8,1) & 105  & 7663 & 42215 & [2, 2, 4, 1, 4, 4]\\
(41,30,8,3) & 123  & 3565 & 61217 & [2, 2, 4, 1, 2, 4]\\
(43,34,8,1) & 137  & 14523 & 78151 & [2, 3, 4, 1, 4, 4]\\
(49,22,15,12) & 147  & 1125 & 92880 & [2, 2, 2, 3, 1, 2]\\
(43,20,14,9) & 149  & 2385 & 94336 & [3, 2, 4, 1, 3, 4]\\
(51,28,13,10) & 153  & 1375 & 100171 & [2, 2, 2, 3, 1, 4]\\
(49,38,8,3) & 155  & 6201 & 102969 & [2, 3, 4, 1, 4, 4]\\
(53,26,12,3) & 159  & 6517 & 108079 & [2, 2, 4, 1, 4, 2]\\
(43,32,10,1) & 161  & 26727 & 112336 & [3, 2, 4, 1, 4, 4]\\
(55,21,18,16) & 165  & 1043 & 117928 & [2, 2, 2, 2, 3, 1]\\
(59,42,9,8) & 177  & 2535 & 138768 & [2, 2, 2, 2, 1, 4]\\
(57,36,10,1) & 181  & 26970 & 149200 & [2, 3, 4, 1, 4, 4]\\
(49,36,10,3) & 183  & 11591 & 151842 & [3, 2, 4, 1, 2, 4]\\
(61,20,18,3) & 183  & 8965 & 152011 & [2, 2, 4, 1, 4, 2]\\
(61,22,16,7) & 183  & 3841 & 152019 & [2, 2, 2, 1, 4, 4]\\
(57,26,14,3) & 185  & 10176 & 154800 & [2, 3, 2, 1, 4, 4]\\
(61,38,10,3) & 193  & 10230 & 169630 & [2, 3, 4, 1, 4, 4]\\
(65,46,11,8) & 195  & 2533 & 172715 & [2, 2, 2, 3, 1, 4]\\
(59,22,18,1) & 199  & 38010 & 180617 & [3, 2, 4, 1, 4, 4]\\
(61,23,22,16) & 199  & 1593 & 180630 & [3, 2, 3, 4, 4, 1]
\end{array}\]
\caption[Table 1]{All reduced weight systems 
$(v_1,v_2,v_3,v_4;d)$ with $d\leq 200$, 
$\frac{d}{2}>v_1\geq v_2\geq v_3\geq v_4$,
and condition $\oooo{(C2)}$, but not condition $(C2)$.}
\label{table 1}
\end{table}

The table does not contain Ivlev's example
$(58,33,24,1;265)$ because $265>200$. 
Ivlev's example has Milnor number 66516.
All examples in Table \ref{table 1} have smaller Milnor 
numbers. 
We checked also with PARI/GP for all 103 reduced 
weight systems with $d\leq 360$ and $\oooo{(C2)}$, 
but not $(C2)$, that  $\rho_{({\bf v};d)}\in\N_0[t]$ holds. 
It is not clear whether there 
are weight systems with larger $d$ which satisfy $\oooo{(C2)}$
and $\rho_{({\bf v};d)}\in\Z[t]-\N_0[t]$.

Arnold \cite{AGV85} associated to each reduced weight
system $(v_1,...,v_n;d)$ which satisfies $\oooo{(C2)}$ 
one {\it type} or several {\it types} in the following way. 
A {\it type} is a conjugacy class with respect to the 
symmetric group
$S_n$ of a (not necessarily bijective) map $\sigma:N\to N$.
Here $\sigma$ is a map which satisfies $v_j|(d-v_{\sigma(j)})$
for $j\in N$. 
For each $j\in N$ such a $\sigma(j)$ exists because of
condition $\oooo{(C2)}$ for sets $J\subset N$ with $|J|=1$. 
The number of types is $3$ for $n=2$, 
$7$ for $n=3$ \cite{AGV85},
$19$ for $n=4$ \cite{YS79}, $47$ for $n=5$ and $128$ for
$n=6$ \cite[Examples 3.2]{HK12}. 
In \cite{YS79} the types for $n=4$ are numbered
$I,II,...,XIX$ 
(this is reproduced in \cite[Example 3.2 (iii)]{HK12}). 
At most 9 of the 19 types allow weight systems
which satisfy $\oooo{(C2)}$, but not $(C2)$, the weight
systems $V,VIII,XI,XII,XIII,XV,XVI,XVII$ and $XIX$.
The other 10 types allow only weight systems of 
Thom-Sebastiani sums of cycle type singularities 
and chain type singularities.
The 23 weight systems in table \ref{table 1}
realize 7 of the 9 types, they do not realize $XIII$ and $XIX$. 
\end{remark}

\begin{remark}\label{t6.9}
(i) K. Saito conjectured \cite[(3.13) and (4.2)]{Sa88} that 
the weight system ${\bf w}=(w_1,...,w_n;1)$ of an isolated 
quasihomogeneous singularity with multiplicity $\geq 3$
satisfies $\psi_{\bf w}(d_{\bf w})>0$. 
Here $d=d_{\bf w}$ for the equivalent reduced integer
weight system $(v_1,...,v_n;d)$, because condition $(C2)$
for the set $J=N$ implies $v_j|d$ for any $j$. 
Counter-examples with $n=4$ 
to this conjecture were given in \cite[Examples 5.5]{HZ19}.
The tables in \cite{HK11} show that for $n=4$ and
$\mu\leq 500$ the only counter-examples are
those in \cite[Examples 5.5]{HK12}.
{\it Problem 5 (b)} asked for an answer to the question whether in the case $n=4$ and $\mu\geq 501$ the only counter-examples
are those in \cite[Examples 5.5]{HZ19}. We did not settle this.

\medskip
(ii) With the software PARI/GP \cite{PARI19} we checked that
for $n=5$ and $d\leq 200$ only 10 counter-examples to the 
conjecture in part (i) exist. They are all obvious 
generalizations of the examples in \cite[Examples 5.5]{HZ19}. 
They are given in the following table.
They are Thom-Sebastiani sums of singularities of
type $A_k$, $D_k$ (with odd $k$), and/or a chain type 
singularity with one of the three normalized weight systems
$$
K_1:\ \bigl(\frac{1}{4},\frac{3}{20},\frac{17}{40}\bigr),\quad 
K_2:\ \bigl(\frac{1}{4},\frac{3}{20},\frac{25}{56}\bigr),\quad 
K_3:\ \bigl(\frac{1}{6},\frac{5}{18},\frac{13}{36}\bigr).$$
(See \cite[ch. 4 and 5]{HZ19} for formulas around the
singularities $A_k$, $D_k$ and the chain type singularities.)
The number $L$ has the same meaning as in Remark \ref{t6.8}.
\begin{table}
\[\begin{array}{l|r|r|r|r}
\textup{Thom-Sebastiani sum} 
& (v_1,v_2,v_3,v_4,v_5) & d  & \mu & L \\
\hline 
D_{13}\otimes K_1
&(55,51,30,18,10) & 120 & 299 & 305495\\
D_{13}\otimes D_{21}\otimes A_3
&(57,55,30,10,6) & 120 & 819 & 307654\\
D_{13}\otimes K_2
&(77,75,42,18,14) & 168 & 403 & 740415\\
D_{11}\otimes D_{29}\otimes A_3
&(81,77,42,14,6) & 168 & 1131 & 744420\\
D_{11}\otimes K_3
&(81,65,50,30,18) & 180 & 253 & 903983\\
D_{11}\otimes D_{19}\otimes A_5
&(85,81,30,18,10) & 180 & 1045 & 907794\\
D_{11}\otimes D_{19}\otimes A_2
&(85,81,60,18,10) & 180 & 418 & 907798\\
D_{31}\otimes K_3
&(87,65,50,30,6) & 180 & 713 & 909833\\
D_{19}\otimes D_{31}\otimes A_5
&(87,85,30,10,6) & 180 & 2945 & 909911\\
D_{19}\otimes D_{31}\otimes A_2
&(87,85,60,10,6) & 180 & 1178 & 909915
\end{array}\]
\caption[Table 2]{All reduced weight systems 
$(v_1,v_2,v_3,v_4,v_5;d)$ with $d\leq 200$, 
$\frac{d}{2}>v_1\geq ...\geq v_5$, condition $\oooo{(C2)}$, 
and $\psi_{\bf w}(d_{\bf w})=0$.}
\label{table 2}
\end{table}
We calculated with PARI/GP also the number of all reduced 
weight systems $(v_1,v_2,v_3,v_4,v_5;d)$ with 
$\oooo{(C2)}$ and $d\leq 200$ and 
$\frac{d}{2}>v_1\geq ...\geq v_5$. It is 1176435. 

\medskip
(iii) K. Saito posed also the weaker conjecture
\cite[(3.13) and (4.2)]{Sa88} that the weight system 
${\bf w}$ of any isolated quasihomogeneous singularity 
satisfies $\psi_{\bf w}(d_{\bf w})>0$ or
$\psi_{\bf w}(d_{\bf w}/2)>0$. 
{\it Problem 5 (a)} in \cite{HZ19} rephrases this conjecture.
If the only counter-examples to the stronger conjecture in
part (i) are the obvious generalizations of the
examples 5.5 in \cite{HZ19}, the weaker conjecture is true. 
But this is not clear.
\end{remark}




\begin{thebibliography}{AAA9}
\bibitem[Ai79]{Ai79} M. Aigner: \quad 
    Combinatorial theory. Grundlehren der math. Wiss. 234,
    Springer Verlag, Berlin Heidelberg New York, 1979.
\bibitem[AGV85]{AGV85} V.I. Arnold, S.M. Gusein-Zade, 
    A.N. Varchenko: \quad
    Singularities of differentiable maps, Volume I.
    Birkh\"auser, Boston, Basel, Stuttgart, 1985.
\bibitem[Fl00]{Fl00} A.R. Fletcher: \quad
    Working with weighted complete intersections.
    In: Explicit birational geometry of 3-folds
    (eds. A. Corti, M. Reid), 
    London Mathematical Society Lecture Note Series 281,
    Cambridge University Press, 2000, pages 101--174.
\bibitem[He11]{He11} C. Hertling: \quad
    $\mu$-constant monodromy groups and marked singularities.
    Ann. Inst. Fourier, Grenoble {\bf 61.7} (2011), 2643--2680.
\bibitem[He20]{He20} C. Hertling: \quad 
    Automorphisms with eigenvalues in $S^1$ of a $\Z$-lattice 
    with cyclic finite monodromy. 
    Acta Arithmetica {\bf 192.1} (2020), 1--30.
\bibitem[HK11]{HK11} C. Hertling, R. Kurbel: \quad
    Tables of weight systems of quasihomogeneous
    singularities. 15.08.2011. 
    {On the homepage: 
    hilbert.math.uni-mannheim.de/CQS-homepage/index.html.}
\bibitem[HK12]{HK12} C. Hertling, R. Kurbel: \quad
    On the classification of quasihomogeneous singularities.
    Journal of singularities {\bf 4} (2012), 131--153.
\bibitem[HM20-1]{HM20-1} C. Hertling, M. Mase: \quad 
    The integral monodromy of cycle type singularities.
    Preprint, arXiv:2009.07533v1, 38 pages, 16.09.2020. 
\bibitem[HM20-2]{HM20-2} C. Hertling, M. Mase: \quad 
    The integral monodromy of isolated quasihomogeneous 
    singularities. Accepted for publication in the
    journal Algebra and Number Theory.
    Preprint version arXiv:2009.08053v1, 84 pages, 
    17.09.2020.
\bibitem[HZ19]{HZ19} C. Hertling, Ph. Zilke: \quad 
    Seven combinatorial problems around isolated 
    quasihomogeneous singularities.
    Journal of Algebraic Combinatorics {\bf 50} (2019),
    447--482.
\bibitem[Ko76]{Ko76} A.G. Kouchnirenko:  \quad
    Poly\`edres de Newton et nombres de Milnor.
    Invent. Math. {\bf 32} (1976), 1--31.
\bibitem[Ko77]{Ko77} A.G. Kouchnirenko:  \quad
    Criteria for the existence of a non-degenerate 
    quasihomogeneous function with given weights. 
    (In Russian.) Usp. Mat. Nauk {\bf 32:3} (1977), 169--170.
\bibitem[KS92]{KS92} M. Kreuzer, H. Skarke: \quad
    On the classification of quasihomogeneous functions.
    Commun. Math. Phys. {\bf 150} (1992), 137--147.
\bibitem[Mi68]{Mi68} J. Milnor: \quad
    Singular points of complex hypersurfaces.
    Ann. of Math. Stud. {\bf 61}, Princeton Univ. Press, 1968.
\bibitem[MO70]{MO70} J. Milnor, P. Orlik: \quad 
    Isolated singularities defined by weighted homogeneous 
    polynomials. Topology {\bf 9} (1970), 385--393.
\bibitem[Or72]{Or72} P. Orlik: \quad
    On the homology of weighted homogeneous manifolds.
    In: Lecture Notes in Math. 298, Springer, Berlin, 1972, 
    260--269.
\bibitem[OR76]{OR76} P. Orlik, R. Randell: \quad
    The classification and monodromy of weighted 
    homogeneous singularities.
    Preprint, 1976 or 1977, 40 pages.
\bibitem[PARI19]{PARI19}
    The PARI~Group, PARI/GP version \texttt{2.11.2}, 
    Univ. Bordeaux, 2019,
    \url{http://pari.math.u-bordeaux.fr/}.
\bibitem[Sa71]{Sa71} K. Saito: \quad
    Quasihomogene isolierte Singularit\"aten von 
    Hyperfl\"achen. Invent. Math. {\bf 14} (1971), 123--142.
\bibitem[Sa87]{Sa87} K. Saito: \quad
    Regular systems of weights and their associated 
    singularities.
    In: Complex analytic singularities. 
    Advanced Studies in Pure Math. {\bf 8}, 
    Kinokuniya \& North Holland, 1987, 479--526.
\bibitem[Sa88]{Sa88} K. Saito: \quad 
    On the existence of exponents prime to the Coxeter number.
    J. Alg. {\bf 114} (1988), 333-356.
\bibitem[ST71]{ST71} M. Sebastiani, R. Thom: \quad 
    Un resultat sur la monodromie.
    Invent. Math. {\bf 13} (1971), 90--96.
\bibitem[Sh79]{Sh79} O.P. Shcherbak: \quad 
    Conditions for the existence of a non-degenerate mapping
    with a given support. Func. Anal. Appl. {\bf 13} (1979), 
    154--155.
\bibitem[YS79]{YS79} E. Yoshinaga, M. Suzuki: \quad
    Normal forms of nondegenerate quasihomogeneous
    functions with inner modality $\leq 4$.
    Invent. Math. {\bf 55} (1979), 185--206. 
\bibitem[Wa96]{Wa96} C.T.C. Wall: \quad 
    Weighted homogeneous complete intersections.
    In: Algebraic geometry and singularities 
    (La R\'abida, 1991). Progr. Math. {\bf 134}, 
    Birkh\"auser, Basel, 1996, 277--300.
\end{thebibliography}
\end{document}